\documentclass[review]{elsarticle}

\usepackage{lineno,hyperref}
\modulolinenumbers[5]

\journal{Journal of \LaTeX\ Templates}

\usepackage{amsfonts,amssymb,amsbsy}
\usepackage{latexsym}
\usepackage{amsmath}
\usepackage{xcolor}

\newtheorem{theorem}{Theorem}

\newtheorem{corollary}{Corollary}
\newtheorem{example}{Example}
\newtheorem{remark}{Remark}
\newtheorem{definition}{Definition}

\newcommand{\R}{\mathbb{R}}

\newcommand{\N}{\mathbb{N}}

\begin{document}

\begin{frontmatter}

\title{Exponential stability analysis for a  class of switched nonlinear time-varying functional differential systems}
\title{Exponential stability analysis for a  class of switched nonlinear time-varying functional differential systems\tnoteref{mytitlenote}}

\author[mymainaddress]{Nguyen Khoa Son\corref{mycorrespondingauthor}}
	\cortext[mycorrespondingauthor]{Corresponding author}
	\address[mymainaddress]{Institute of Mathematics, Vietnam Academy of Science and Technology, 18 Hoang Quoc Viet Rd., Hanoi, Vietnamm}\ead{nkson@vast.vn}
	
	\author[mysecondaryaddress]{Le Van Ngoc}
		
	\address[mysecondaryaddress]{Department of Scientific Fundamentals, Posts and Telecommunications Institute of Technology, Km10 Nguyen Trai Rd, Hanoi, Vietnam}
	\ead{ngoclv@ptit.edu.vn}


\begin{abstract}
This paper proposes a unified  approach for  studying global exponential stability of  a general class of switched  systems described by time-varying nonlinear functional differential equations. Some new delay-independent criteria of global exponential stability are  established for this class of systems under arbitrary switching which satisfies some assumptions on  the average dwell time. The obtained criteria are shown to cover and improve many previously known results, including, in particular, sufficient conditions for absolute exponential  stability of switched time-delay systems with sector nonlinearities. Some simple examples are given to illustrate the proposed method.
\end{abstract}
\end{frontmatter}

\section{Introduction}\label{Section1}

Functional differential equations (or FDEs, for short) and, in particular, time-delay differential equations  have numerous applications in science and engineering, see, e.g., \cite{Hale}. In the framework of qualitative theory, stability of FDEs is always an important research topic that has been studied intensively during the past decades, see e.g.  \cite{Mi_Ni, ngoc_tinh}  and the references therein. Recently, stability problems have been  considered intensively also for time-delay switched systems. Recall that a switched system is a type of hybrid dynamic system which consists of a family of subsystems and a rule called a switching signal that chooses an active subsystem from the family at every instant of time. The reader is referred to \cite{Liberzon} and also \cite{Shorten, Lin}  for more details on stability problems of switched systems. Some recent development in exponential stability analysis of time-delay switched systems can be found, for instance, in \cite{Kim, Wang, Li2017, Li2018,  Liu2018, Gao_Lib, Tian_Sun, sun_wang2013}, where the Lyapunov-Krasovskii functionals method and the  comparison principle have been widely used. Similar problems have been investigated also for different classes of switched positive systems, see e.g. \cite{Dong2015, mengIEEE16, zhao_IEEE17, wang2020, alex_mason2014}.

The main purpose of this paper is to propose an unified  approach for studying exponential stability of a  general class of time-varying switched systems, described by nonlinear functional differential equations, that is based on the comparison principle and the average dwell time (ADT) switching concept \cite{Hes_Mor, Liberzon}. As the main results, we will establish some verifiable criteria for exponential stability of the zero solution of switched nonlinear FDE system under arbitrary switching, satisfying some ADT assumptions. These results, to the best of our  knowledge, are new in the literature and, moreover, cover many previous results on exponential stability of time-delay switched systems. In particular, these results apply to switched time-delay systems with sector nonlinearities, establishing some new sufficient conditions for absolute exponential stability of the zero solution, that relax and extend the conditions given
in \cite{sun_wang2013, alex_mason2014, alex_zab2017, zhang_zhao2016}.  

The layout of the paper is  as follows.  In Section 2, we present our main results for  a general class of switched time-varying systems that is modeled by nonlinear FDE, using the comparison principle and ADT switching laws.  In Section 3, we show the application of the main results to classes of switched linear  systems with time-delay, while the absolute exponential stability of systems with sector nonlinearities is considered in Section 4. Two numerical examples are given in Section 5 to illustrate the use of the obtained results and a conclusion is given in Section 6 to summarize the contribution of the paper and mention about some possible further extensions of the results. 

The following notation will be used throughout the paper. $\R, \R_+$ and $\N_+$ will stand for the sets of real numbers, non-negative real numbers and non-negative integers, respectively. For  $r\in \N_+, \underline r$ denotes the set of numbers $ \{1,2,\ldots, r \}$. For vectors $x=(x_i), y=(y_i) \in \R^n$  we write $x\geq y$ and $x\gg y$ iff $x_i \geq y_i$ and $x_i> y_i$, for all $i\in \underline n,$ respectively. Denote $|x|=(|x_i|)$ and $x^{\top}$ is the transpose of $x$. Similar notation is applied for matrices, that is, for $A=(a_{ij}), B=(b_{ij})\in \R^{n\times n}$ we write, for instance,  $A\gg B$ iff $a_{ij}>b{ij}, \forall i,j \in \underline n$ and $|A|=(|a_{ij}|)$. Without loss of generality, the norm of vectors $x\in \R^n$ is assumed to be the $\infty$-norm: $$\|x\|=\|x\|_{\infty}=\max_{1\leq i\leq n}|x_i|.$$   
For $h>0, \mathcal{C} := C([-h, 0], \R^n)$ denotes the Banach space of continuous functions $\varphi: [-h, 0] \rightarrow \R^n $  with the norm $\|\varphi\|= \max_{\theta\in[-h,0]} \|\varphi(\theta)\|$ and $ \mathcal{C}_+:= \{\varphi \in \mathcal{C}: \varphi(\theta)\geq 0, \ \forall \theta\in [-h,0]\}$. 
$NBV([-h,0],\R)$ stands for the linear space of  all normalized functions $\psi: [-h,0] \rightarrow \R$ with bounded variation $ Var([-h,0],\psi)$, so that $\psi$  is left-side continuous on the interval $(-h,0]$ and $\psi (-h) = 0$ and $Var([-h,0],\psi):= \sup_{P[-h,0]}\sum_{k}|\psi(\theta_k)-\psi(\theta_{k-1})| <+\infty$, where the supremum is taken over the set of all finite partitions $P$ of the interval $[-h, 0]$.  It is well-known that, for any $\psi \in NBV([-h,0],\R)$ and any continuous function $\beta \in   C([-h, 0], \R),$ we have
\begin{equation}\label{stielj}
	\int_{-h}^0 d[\psi(\theta)] \beta(\theta) \leq Var([-h,0],\psi) \max_{\theta\in [-h,0]}|\beta(\theta)|,
\end{equation}
where the integral is understood in the sense of Riemann-Stieltjes. Similarly,  $NBV([-h,0],\R^{n\times n})$ will stand for the linear space of all matrix functions $\eta \!:\! [-h,0] \!\rightarrow \!\R^{n\times n}$ such that  $\eta_{ij}(\cdot) \in NBV([-h,0],\R), \forall i,j\in \!\underline \!\underline n.$  Thus, to each $\eta \in NBV([-h,0],\R^{n\times n}) $ we can associate a nonnegative $(n\times n )$-matrix of variations 
\begin{equation}
	\label{var}
	V(\eta):=\big(Var([-h,0], \eta_{k,ij}) \big) \geq 0 .
\end{equation}
\noindent  If the function $f:\R_+\times \R^n\rightarrow \R^n$ is differentiable w.r.t.  $x$  then $ J_xf=(\frac{\partial f_i(t,x)}{\partial x_j})\in \R^{n\times n}$ denotes the Jacobian of $f$. Finally, for any  matrix $A\in \R^{n\times n}$ we associate the {\it Metzler matrix} $\mathcal{M}(A)= (\bar{a}_{ij})$, by setting, for all $i,j \in \underline  n, \bar{a}_{ii}= a_{ii}, $ and $\bar{a}_{ij}= |a_{ij}|$ if $i\not= j$.   
\vskip-1.0cm
\section{Main results}

Consider a  switched system described by nonlinear time-varying functional  differential equations of the form
\begin{equation}
\label{DSwFDS}
\dot{x}(t) = f_{\sigma(t)}(t, x,x_t),\  t\geq 0, \ \sigma\in \Sigma_+,
\end{equation}
where, for each $t\geq 0,\  x_t(\cdot) \in \mathcal C$ is defined by $x_t(\theta) := x(t + \theta ), \theta \in [-h, 0] $  for given
$h > 0$, $f_{\sigma(t)}(\cdot,\cdot,\cdot)\in \mathcal F := \{f_k(\cdot,\cdot,\cdot),  k\in \underline{N}\} $ -a given family of $N$ continuous  nonlinear  maps $f_k(t,x,\varphi):\R_+\times \R^n \times \mathcal C \rightarrow \R^n$ satisfying some assumptions specified below and $\Sigma_+$ is a given set of {\it switching signals}. Each switching signal is assumed to be a piece-wise constant function $\sigma : [0,\infty) \rightarrow \underline {N}$  which is right-side continuous, having points of discontinuity $\tau_k, k=1,2,\ldots,  $ known as the {\it switching instances}. As a common convention in the theory of switched systems, switching signal $\sigma$ is assumed to have on each bounded interval only a finite number discontinuities $\tau_k.$
Such a set $\Sigma_+$ excludes, for instance, any switching signal whose discontinuities have a finite accumulation point (i.e. a signal having Zeno's behaviour) or occur at $\tau_{2j}=j, \tau_{2j+1}=j+\frac{1}{2j+1}, j=0,1,2, \ldots$ for which, clearly,  $\tau_{min}(\sigma)=0$.

\noindent  Thus, each signal $\sigma\in \Sigma_+$ performs  switchings  between the following $N$ nonlinear time-varying functional differential {\it constituent} subsystems  of the form
\begin{equation}
	\label{DSwFDSk}
	\dot{x}(t)= f_k(t,x,x_t), \  t\geq 0,	\ k\in \underline N,
\end{equation}
\noindent Throughout the paper we will assume that the family $\mathcal F$  satisfies the following condition:

\noindent (A1)  For each $k\in \underline  N, f_k(t, 0,0)\!=\!0, \forall t\!\geq\! 0$ and, for any $ \varphi\in \mathcal C$,  the subsystem \eqref{DSwFDSk} admits a global unique solution $x(t)\!=\!x(t,\varphi), t\!\geq\! -h$ satisfying the initial condition $x(\theta)\! =\!\varphi(\theta), \theta\!\in\! [-h,0]$.   

\noindent It is well-known \cite{Hale} that the second assumption in (A1) holds, for instance,  if for each $k\in \underline N, \ f_k (t,x,\varphi) $ is uniformly bounded and Lipschitz continuous in $x,\varphi$ on any compact subsets in $ \R_+\times \R^n\times \mathcal C $.   Under the assumption (A1),  $x(t) \equiv 0, t\geq -h $ is the {\it zero solution} of \eqref{DSwFDS} and, for each $\varphi \in \mathcal{C}$ and each $\sigma \in \Sigma_+$, the system \eqref{DSwFDS} admits a unique global solution $ x(t)\!=\! x(t,\varphi, \sigma ), t\!\geq\! -h $, satisfying the initial condition 
\begin{equation}\label{init}
 x(\theta)= \varphi(\theta),\  \theta\in [-h,0].
 \end{equation}
   Note that the solution $x(t) $ of \eqref{DSwFDS} is a absolutely continuous function on $[0,+\infty)$ and differentiable everywhere, except for the set of switching instances $\{\tau_k\}$ of $\sigma $  where $x(t)$ has only Dini right- and left-derivatives $D^+x(\tau_k), D^-x(\tau_k) $  which are generally different. 
   
  \begin{definition}\label{GES}
  	Let a set of switching signals $\Sigma \subset \Sigma_+$ be given. The zero solution of the switched system \eqref{DSwFDS} is said to be globally exponentially stable (shortly, GES) over $\Sigma $ if there exist positive numbers $M, \alpha $ such that for any $\varphi\in C([-h,0],\R^n) $ and any $\sigma \in \Sigma $ the solutions $ x(t,\varphi, \sigma) $  of \eqref{DSwFDS} satisfies
  	\begin{equation}
  		\label{conditiondef1}
  		\|x(t,\varphi,\sigma)\|\leq M e^{-\alpha t}\|\varphi\|, \quad \forall t\geq 0.
  	\end{equation}
  \end{definition}
	\noindent 
	Obviously, for each $k\in \underline N$,  the  switching law $\sigma (t)\equiv k, t\geq 0,$ belongs to $\Sigma_+$. Therefore, if the zero solution of switched system \eqref{DSwFDS} is GES over the set of switching signals $ \Sigma_+$ then the same property holds true for all constituent subsystems \eqref{DSwFDSk}. The converse, as is well-known,  is generally not true, even for switched linear systems with no delays (see, e.g. \cite{Liberzon}). 
	
	\begin{remark}
{\rm The concept 'GES over a set  of switching signals' $\Sigma$ in Definition \ref{GES} has been introduced and studied firstly in \cite{Hes_Mor}. In many subsequent works, the equivalent concept 'GES under arbitrary switchings' $\sigma\in \Sigma $  is used more frequently. }
	\end{remark}
	
\noindent In the sequels we are also interested in exponential stability of switched positive systems which concept is given by the following 

\begin{definition} \label{positive} The switched nonlinear system of the form \eqref{DSwFDS} is said to be positive if for any initial function $\varphi \in \mathcal{C}_+$ and any switching signal $\sigma \in \Sigma$ the associate solution $x(t)= x(t,\varphi,\sigma)$ satisfies $x(t)\geq 0, \ \forall t\geq 0$. 
\end{definition}

In order to analyze  exponential stability  of  switched nonlinear FDE of the form \eqref{DSwFDS} by using the comparison principle, besides the above standard condition (A1),  we will assume that the following specific conditions hold for any  $f_k\in \mathcal {F}, \ k\in \underline N$.

(A2) $ f_k(t,x,\varphi)$  is continuously differentiable w.r.t. the variable $x$, and there exist  continuous  matrix functions $A_k(t):= \big(a_{k,ij}(t)\big)\in \R^{n\times n}$ such that, for any $t \geq \!0,\ x \in \R^n,\ \varphi \in \mathcal C $,
\begin{equation}\label{A2}
\mathcal M \big(J_xf_k(t,x,\varphi)\big) \leq A_k(t), \ \forall k\in \underline N,
\end{equation}
 or,  equivalently, for all $ \ k\in \underline N, \ i\in \underline n$,  
\begin{equation}\label{A2-1} 
\frac{\partial f_{k,i}}{\partial x_i}(t,x,\varphi)\!\leq\! a_{k,ii}(t);\ \big| \frac{\partial f_{k,i}}{\partial x_j}(t,x,\varphi)\big|\! \leq a_{k,ij}(t),\ \forall j \not= i, \ j \in \underline n.
\end{equation}

(A3) There exists a continuous function $	L_k(\cdot, \cdot): \R_+\times \mathcal C \ \rightarrow\  \R^n$, defined by 
\begin{equation}\label{Lk} 
L_k(t,\varphi)=\int_{-h}^0d_{\theta}[\eta_k(t,\theta)]\varphi(\theta) 
\end{equation}
 such that, for any $t \geq 0, \ \varphi \in \mathcal C $, 
\begin{equation}\label{A3} 
|f_k(t,0,\varphi)| \leq |L_k(t,\varphi)|,\   \forall k\in \underline N, 	
	\end{equation}
where, for each $k\in \underline N$ and each $t\geq 0$, $\eta_k(t,\cdot):=\big(\eta_{k,ij}(t,\cdot)\big)\in NBV([-h,0],\R^{n\times n})$ and the integral \eqref{Lk} is understood in the sense of Riemann-Stieltjes.

\begin{remark}\label{riesz} {\rm
	By the Riesz representation theorem, any continuous function $	L_k(\cdot, \cdot): \R_+\times \mathcal C \ \rightarrow\  \R^n$, which, for each fixed $t$,  is a linear bounded operator with respect to $\varphi \in \mathcal{C}$ can be represented in the form \eqref{Lk} with an appropriate function $\eta_k(t,\cdot)\in NBV([-h,0],\R^{n\times n})$. Thus, Assumption (A3) amounts to saying that nonlinear functions $f_k(t,0,\cdot)$ are overbounded by some linear operators $L_k(t,\cdot)$  for each $t\geq 0$.}
\end{remark}

We recall (see, e.g. \cite{Hes_Mor}) that, for given two numbers $\tau_a >0, N_0\geq 0$, a switching signal $\sigma \in \Sigma_+ $ is said to have an {\it average dwell time} (or ADT, for short)  $\tau_a $, with the chatter bound $N_0$, if for any $t >0$ the number $N_{\sigma}(0,t)$ of discontinuities of $\sigma$ on the interval  $(0,t]$ satisfies
\begin{equation} \label{ADT}
	N_{\sigma}(0,t) \leq N_0+ \frac{t}{\tau_a}.
\end{equation}
The set of all switching signals having ADT $\tau_a $ is denoted by $\Sigma_{\tau_a,N_0}$ or simply by $\Sigma_{\tau_a}$ if $N_0=0$.  It follows that for any $\sigma \in \Sigma_{\tau_a,N_0}$, if we ignore the first $N_0$ discontinuities, then   the average dwell time between any two consecutive switching instances  is at least $\tau_a. $ Obviously, for each fixed $N_0\geq 0$, if  $\tau_1>\tau_2 >0$, then
$$
\Sigma_{\tau_1,N_0}\subset \Sigma_{\tau_2,N_0}\subset  \Sigma_{+}. 
$$
Therefore, if the system \eqref{DSwFDS} is GES over  $ \Sigma_{+}$ it is also GES over $\Sigma_{\tau_a,N_0}$ for any $\tau_a >0$. On the other hand, it is a well known fact in the theory of switched systems  that if each
subsystem is GES and the average dwell
time $\tau_a$ of the switching signal is sufficiently large, then the switched
system is GES. This raised a problem of remarkable interest to find a lower bound for ADT $\tau_a$ which guarantees GES for the switched system over the set of switching signals $\Sigma_{\tau_a,N_0}$. The earliest result in this line of research has been derived in \cite{Morse, Hes_Mor}, for linear systems of the form $\dot x(t) = A_{\sigma(t)}x(t)$. Subsequently, this problem has been considered for other types of switched systems, including linear time-delay systems and nonlinear systems, based  either on  the Lyapunov functions method, like in the last mentioned works, or on the comparison principle (see, e.g. \cite{Gao_Lib}, \cite{Tian_Sun} and the references therein).

\noindent We are now in a position to prove the main result of this paper which gives a lower bound for such ADT of switching signals, for the general class of nonlinear switched FDE of the form \eqref{DSwFDS}. 
\begin{theorem}\label{main}
	Consider the switched nonlinear  FDE system \eqref{DSwFDS} which satisfies (A1,A2,A3). Assume, moreover, that there exist vectors $\xi_k\in \R^n, \xi_k \gg 0, k\in \underline N$ and a positive number $\alpha > 0$ such that
	\begin{equation} 
\label{cond7}
\big(A_k(t)+ e^{\alpha h} V(\eta_k(t,\cdot))\big)\xi_k \ll -\alpha \xi_k, \ \forall t\geq 0,\  \forall k \in \underline N. 
\end{equation} 
Then, for each fixed chatter bound $N_0\geq 0$,   the zero solution of \eqref{DSwFDS} is GES over the set   $ \Sigma_{\tau_a,N_0}$ of switching signals with ADT $\tau_a$ satisfying
\begin{equation}\label{tau}
\tau _a > \tau_*:=\frac{\ln \; \gamma}{\alpha},
\end{equation}
where 
\begin{equation}\label{mu}
	\gamma := \max\big\{\frac{\xi_{k,i}}{\xi_{l,i}}: k,l\in \underline N, i\in \underline n\big\},  \ \xi_k := (\xi_{k,1}\  \xi_{k,2}\ \cdots\ \xi_{k,n})^{\top}.  
\end{equation}
\end{theorem}
\vskip -0.3cm
{\it Proof. } 
Without loss of generality we can assume that $ \|\xi_k\|=1,  k \in \underline N$ for vectors $\xi_k \gg 0 $ in \eqref{cond7}.  Let $\tau_a$ be chosen to satisfy \eqref{tau} and $\sigma \in \Sigma_{\tau_a}$ be an arbitrary switching signal with switching instances $0=\tau_0< \tau_1< \ldots < \tau_p <\tau_{p+1} <\ldots$. Thus, for any $t>0$ the number of discontinuities of $\sigma$ on the inverval $(0,t]$  is $ N_{\sigma}(0,t) = p$,  whenever $t\in [\tau_p,\tau_{p+1}).$  For an arbitrary  $\varphi \in \mathcal C$, let $x(t)=x(t,\varphi,\sigma) $ be the corresponding solution of \eqref{DSwFDS} satisfying the initial condition \eqref{init}. Assume that $\sigma(\tau_p)= l_p\in \underline N,$ i.e. the subsystem  $\dot x(t)= f_{l_p}(t,x,x_t)$ is active on $[\tau_p,\tau_{p+1}), p=0,1,\ldots$. For any $\delta>1$, define the piece-wise continuous functions $y_i(t), t\geq -h, i\in \underline n$, by setting 
\begin{equation}\label{y}	
	y_i(t)= 
	\begin{cases}
	 M_{\delta}e^{-\alpha t}\xi_{l_0,i}\|\varphi\|  \ &\text{if}\   t\in [-h,\tau_0],\\
	M_{\delta}e^{-\alpha t}\xi_{l_p,i}\|\varphi\| &\text{if} \ t \in [\tau_p, \tau_{p+1}), \ p=0,1,2,\ldots
	\end{cases}
\end{equation}
where $M_{\delta}=\delta.\gamma$. Since, obviously,  $M_{\delta}>1 $, it follows from  \eqref{init} and \eqref{y} that,  for each $i\in \underline n$,
\begin{equation}	\label{defk}
	|x_{i}(t)|=|\varphi_{i}(t)|  < M_{\delta}e^{-\alpha t}\xi_{l_0,i}\|\varphi\|=y_{i}(t), \  \forall t\in [-h,0]. 
\end{equation} 
We will show that, for each $ i\in \underline n,  $ 
\begin{equation}
\label{assu1}
|x_i(t)|\leq \gamma ^p y_i(t)=\gamma^pM_{\delta}e^{-\alpha t}\xi_{l_p,i}\|\varphi\| , \  \forall t\in [\tau_p,\tau_{p+1}),\  p=0,1,\ldots
\end{equation}
The proof will be proceeded inductively  on each interval $[\tau_p,\tau_{p+1})$ where the $l_p$-th subsystem is active.   Assume to the contrary that \eqref{assu1} does not hold for $p=0$, then, due to \eqref{defk}, by the continuity, there exist $ i_0\in \underline n, \   \bar t_0\in [\tau_0,\tau_1)$ and $ \varepsilon >0$ such that    
\begin{equation}\label{conditi_The1}
|x_{i}(t)|< y_{i}(t), \forall t\in [-h,\bar t_0), \ \forall i\in \underline n,
\end{equation}
and
\begin{equation}
\label{conditi_Theo2}
|x_{i_0}(\bar t_0)|=y_{i_0}(\bar t_0),
|x_{i_0}(t)| >y_{i_0}(t),\ \forall t\in (\bar t_0,\bar t_0+ \varepsilon).
\end{equation}
Since $\sigma \in \Sigma_{\tau_a}\subset \Sigma_{+}$ and $\bar t_0\in [\tau_0,\; \tau_1)$, the number $ \varepsilon$ in \eqref{conditi_Theo2} can be, moreover, chosen sufficiently small so that $  [\bar t_0,\ \bar t_0+ \varepsilon) \subset [\tau_0,\; \tau_1)$.  Since $\sigma(\tau_0)=l_0$, it follows that $x(t)$ satisfies the subsystem $\dot x(t)= f_{l_0}(t,x,x_t)$ on the interval $ [\tau_0,\tau_1)$.  Denote elements of $A_{l_0}(t)$ and $\eta_{l_0}(t,\cdot)$, respectively, by $ a_{l_0, ij}(t)$ and $ \eta_{l_0,ij}(t,\cdot), i,j \in \underline n$. Then, by the assumption (A2) and the mean value theorem, we have for any  $t\in [\tau_0,\tau_1)$ and  every $i\in\underline{n}$,
\begin{align}\label{mean0}
\dot x_i(t)&= 	f_{l_0,i}(t,x(t),x_t)= f_{l_0,i}(t,x(t),x_t)-f_{l_0,i}(t,0,x_t) + f_{l_0,i}(t,0,x_t)\notag\\
&=\sum_{j=1}^n\bigg(\int_{0}^1 \frac{\partial f_{l_0,i}}{\partial x_j}(t,sx(t),x_t)ds\bigg)x_j(t) + f_{l_0,i}(t,0,x_t).
\end{align}
Therefore, by using \eqref{stielj}, \eqref{A2-1},\eqref{A3}, \eqref{mean0}, we can deduce (for any  $t\in [\tau_0,\tau_1)$ and  every $i\in\underline{n}$ )
\begin{align}
D^+|x_i(t)| & \leq \text{sign}(x_i(t))\dot x_i(t)\notag \\
&\stackrel{\eqref{A2-1},\eqref{mean0}} \leq a_{l_0,ii}(t)|x_i(t)|+\!\sum_{j=1,j\neq i}^n\!a_{l_0,ij}(t)\;|x_j(t)| + |f_{l_0,i}(t,0,x_t)|  \notag\\
& \stackrel{\eqref{A3}}\leq \sum_{j=1}^n a_{l_0,ij}(t)\;|x_j(t)| + \sum_{j=1}^n \bigg| \int_{-h}^0 d_{\theta}[\eta_{l_0,ij}(t,\theta)]x_j(t\!+\!\theta)\bigg| \notag\\
&\stackrel{\eqref{stielj}}\leq \sum_{j=1}^n a_{l_0,ij}(t)\;|x_j(t)|+ \sum_{j=1}^nV(\eta_{l_0,ij}	(t,\cdot))\max_{\theta\in [-h,0]}|x_j(t+\theta)|.
\label{concl1}
\end{align}
It follows easily from \eqref{conditi_The1} that $\max_{\theta\in [-h,0]}|x_j( \bar t_0+\theta)|\leq M_{\delta} e^{-\alpha (\bar t_0-h)}\xi_{l_0,j}\|\varphi\|, \ \forall j\in \underline n$.
Therefore, from \eqref{concl1} (with  $i=i_0, t=\bar t_0) $ and \eqref{cond7} (with $k=l_0$) and using the equality in \eqref{conditi_Theo2}, we get
\begin{equation*}
D^+\!|x_{i_0}(\bar t_0)|\leq M_{\delta} e^{-\alpha \bar t_0}\|\varphi\| \bigg(A_{l_0}(\bar t_0)\xi_{l_0} + e^{\alpha h}V({\eta}_{l_0}(\bar t_0,\cdot)\xi_{l_0}\bigg)_{i_0} 
<M_{\delta}e^{-\alpha\bar t_0}\|\varphi\|(-\alpha\xi_{l_0})_{i_0}\!=\dfrac{d}{dt}y_{i_0}(\bar t_0).
\end{equation*}
On the other hand, by definition of  the Dini right-derivative and the inequality in \eqref{conditi_Theo2}, we have
\begin{align*}
D^+|x_{i_0}(\bar t_0)|= \lim_{\beta\rightarrow0^+}\dfrac{|x_{i_0}(\bar t_0+\beta)|-|x_{i_0}(\bar t_0)|}{\beta}	\geq \lim_{\beta\rightarrow0^+}\dfrac{y_{i_0}(\bar t_0+\beta)-y_{i_0}(\bar t_0)}{\beta}=\dfrac{d}{dt}y_{i_0}(\bar t_0), 
\end{align*}
a contradiction. This  proves that \eqref{assu1}holds for $p=0$. 
Letting $t$ tend to $\tau_1$ and $\delta$ tend to 1 we get from \eqref{mu}, \eqref{y}  that, for each $ i\in \underline n$, 
\begin{align}\label{mu2}
|x_i(\tau_1)|\!\leq \!M_{1}e^{-\alpha \tau_1}\xi_{l_0,i}\|\varphi\|\!=\! \frac{\xi_{l_0,i}}{\xi_{l_1,i}} M_{1}e^{-\alpha \tau_1}\xi_{l_1,i}\|\varphi\|\!< \gamma M_{\delta} e^{-\alpha \tau_1}\xi_{l_1,i}\|\varphi\|\!=\!\gamma\; y_i(\tau_1).
\end{align}
Further, using the  above strict inequality, it can be similarly proved that \eqref{assu1} holds for $p=1$, i.e. 
\begin{equation}\label{tau_3}
|x_i(t)|\leq \gamma y_i(t)= \gamma M_{\delta}e^{-\alpha t}\xi_{l_1,i} \|\varphi\|, \ \forall t\in [\tau_1,\tau_2), \ \forall i\in \underline n.
\end{equation}
Although the proof is mostly repeated as the previous step for \eqref{assu1}, we give it here in  details, for the reader's convenience.   Assume again that \eqref{tau_3} does not hold then,  by continuity,   there exist $ i_1\in \underline n, \   \bar t_1\in (\tau_1,\tau_2)$ and $ \epsilon >0$ such that    
\begin{equation}\label{conditionproof3}
	|x_{i}(t)|< \gamma y_{i}(t), \forall t\in [\tau_1,\bar t_1), \ \forall i\in \underline n,
\end{equation}
and 
\begin{equation}
	\label{conditionproof4}
	|x_{i_1}(\bar t_1)|=\gamma y_{i_1}(\bar t_1),
	|x_{i_1}(t)| > \gamma y_{i_1}(t),\ \forall t\in (\bar t_1,\bar t_0+\epsilon).
\end{equation}
Noticing that the $l_1$-th subsystem  is active on $[\tau_1,\tau_2) $, similarly as the above step, we can estimate the Dini right-derivative $D_+|x(\bar t_1)|$ as (compare with \eqref{concl1})
\begin{equation}\label{D++}
	D^+|x_{i_1}(\bar t_1)| \leq  \sum_{j=1}^na_{l_1,i_1j}(\bar t_1)\;|x_j(\bar t_1)| + \sum_{j=1}^n V({\eta}_{l_1,i_1j})\max_{\theta\in[-h,0]}|x_j(\bar t_1\!+\theta)|.
\end{equation}
Here,  we have to consider two cases: $\tau_1 \leq \bar t_1-h $ and $\bar t_1-h <\tau_1$.  In the first case, it follows from \eqref{conditionproof3} that, for any $\theta \in [-h,0], \ |x_j(\bar t_1+\theta)|\leq  \gamma y_j(\bar t_1 + \theta) \leq \gamma y_j(\bar t_1
-h)  = \gamma M_{\delta} e^{-\alpha (\bar t_1-h)}\xi_{l_1,j}, \ \forall j\in \underline n$. In the second case, by using \eqref{assu1} (with $p=0$) and \eqref{mu}, we have, for any $\theta \in [-h,0]$, 
$$
|x_j(\bar t_1 +\theta)| \leq y_j(\bar t_1 +\theta)= M_{\delta}e^{-\alpha(\bar t_1+\theta)}\xi_{l_0,j} \leq M_{\delta}e^{-\alpha(\bar t_1-h)}\xi_{l_0,j}\leq  \gamma M_{\delta} e^{-\alpha (\bar t_1-h)}\xi_{l_1,j}, \ \forall j\in \underline n.
$$  
Thus in both case, using the equality in \eqref{conditionproof4} and \eqref{D++}, \eqref{cond7} we get the following estimate 
\begin{equation}
	D^+|x_{i_1}(\bar t_1)|\leq \gamma M_{\delta} e^{-\alpha \bar t_1}\bigg(A_{l_1}(\bar t_1)\xi_{l_1}+e^{\alpha h}V({\eta}_{l_1}(\bar t_1,\cdot))\xi_{l_1}\bigg)_{i_1} <\gamma M_{\delta}e^{-\alpha \bar t_1}(-\alpha\xi_{i_1})\!= \gamma \dfrac{d}{dt}y_{i_1}(\bar t_1).
\end{equation}
On the other hand, by the inequality in \eqref{conditionproof4} it follows readily  that  $D^+\!|x_{i_1}(\bar t_1)|  \geq \gamma \dfrac{d}{dt}y_{i_1}(\bar t_1)$, a contradiction. Thus, \eqref{assu1} is proved for $p=1$, that is \eqref{tau_3} holds.

 Again, by letting  $t\rightarrow \tau_2, \delta\rightarrow 1$ in \eqref{tau_3} we get $|x_i(\tau_2))|< \gamma ^2\; M_{\delta} e^{-\alpha (\tau_2-\tau_1)}\xi_{l_1,i} = \gamma^2\; y_i(\tau_2),\ \forall i\in \underline n, $
which similarly implies $|x_{i}(t)|\leq \gamma^2 y_{i}(t), \forall t\in [\tau_2,\tau_3), \ \forall i\in \underline n.$ 
Proceeding as above steps, by induction,  we conclude that \eqref{assu1} holds, which implies immediately
\begin{equation}\label{p}
|x_i(t)| \leq \gamma^p y_i(t)= \gamma^p M_\delta e^{-\alpha t}\xi_{l_p,i}\|\varphi\| \leq \gamma^p M_\delta e^{-\alpha t}\|\varphi\| , \ t \in [\tau_p,\tau_{p+1}), \ i\in \underline n, \  p=0,1,2, \ldots 
\end{equation} 
taking into account that  $\|\xi_k\|=\max_{i\in \underline n}|\xi_{k,i}|=1, \forall k \in \underline N$. Therefore, by  the assumption that $\sigma \in \Sigma_{\tau_a,N_0}$ with ADT $\tau_a$ satifying \eqref{tau}, we can deduce from  \eqref{p} that, for each $i\in \underline n$, 
\begin{equation}
|x_i(t)| \leq \gamma^{N_\sigma(t,0)} M_\delta e^{-\alpha t}\|\varphi\| \leq M_{\delta}e^{N_\sigma(t,0)\ln \gamma}e^{-\alpha t}\|\varphi\|\leq L_\delta e^{(\frac{\ln \gamma}{\tau_a}-\alpha)t}\|\varphi\|,\ \forall t >0,
\end{equation}
where  $L_{\delta}:= M_{\delta}e^{N_0\ln\gamma}>1 $ and $\frac{\ln \gamma}{\tau_a}-\alpha <0  $. This completes the proof.\hspace{6.5cm} $\square$

It is worth noting  that if vectors $\xi_k $ in \eqref{cond7} can be chosen identical, e.g. $\xi_k=\xi \gg 0, \forall k\in \underline N$  then $\gamma =1$ and so $\tau_*=0$, in view of \eqref{tau}.  Hence, by Theorem \ref{main}, the switched system \eqref{DSwFDS} is GES over the set $ \Sigma_{\tau_a,N_0}$   for any ADT $\tau_a >0$ and thus over $\Sigma_+$. Therefore, we have the following consequence of Theorem \ref{main} which has its own interest.

\begin{theorem}\label{main1}
	Consider the switched nonlinear  FDE system \eqref{DSwFDS} which satisfies (A1,A2,A3). Assume, moreover, that there exist vectors $\xi\in \R^n, \xi \gg 0 $  and a positive number $\alpha > 0$ such that
	\begin{equation} 
		\label{cond8}
		\big(A_k(t)+ e^{\alpha h} V(\eta_k(t,\cdot))\big)\xi \ll -\alpha \xi, \ \forall t\geq 0,\  \forall k \in \underline N. 
	\end{equation} 
	Then  the zero solution of \eqref{DSwFDS} is GES over the set   of switching signals	 $ \Sigma_+$.
\end{theorem} 

The Assumption (A2) may be too restrictive to be used in some applications where the functions  $f_k(t,x,\varphi)$ are not necessarily  differentiable. To relax this conservatism we consider the case when the switched system is described by the separable FDE of the form
\begin{equation}\label{separable}
	\dot x(t)= f_{\sigma(t)}(t,x)+ g_{\sigma(t)}(t,x_t), \ t\geq 0,	
\end{equation}
where $f_k:\R_+\times \R^n \rightarrow \R^n,\  g_k: \R_+ \times \mathcal{C} \rightarrow  \R^n, k\in\underline N,$ are continuous functions such that, for each $k\in \underline N, f_k(t,0)=g_k(t,0)=0, \forall t\geq 0 $. Instead of (A2), we will assume that the functions $f_k$ satisfy the folowing assumption which is closely related to the Lipschitz condition.

(A2*) There exist  continuous  matrix functions $A_k(t)\!:=\! \big(a_{k,ij}(t)\big)\!\in\! \R^{n\times n}$ such that, for each $t\!\geq\! 0,\ k\in \underline N,\  A_k(t)$ is a Metzler matrix, i.e. 
\begin{equation}\label{metzA}
\mathcal{M}(A_k(t)) = A_k(t),\ \forall t\geq 0,\  \forall k\in \underline N,  	
\end{equation}
and, for each $ i\in \underline n,\ k\in \underline N$,  
\begin{equation}\label{A2*}
f_{k,i}(t,x)\  \text{sign}\;x_i \leq \sum_{j=1}^n  a_{k,ij}(t)|x_j|,\ \text{for all}\ t\geq 0, x \in \R^n.
\end{equation}
 
\begin{theorem}\label{main2}
	Consider the switched nonlinear  FDE system \eqref{separable} where $f_k(t,x)$  satisfy (A2*) and $g_k(t,\varphi)$ satisfy (A3) i.e. there exist continuous functions of linear operators $L_k(t,\cdot): \mathcal {C} \rightarrow \R^n$ such that 
\begin{equation}\label{gk}
	|g_k(t,\varphi)| \leq |L_k(t,\varphi)|, \ \forall t \geq 0, \ \forall\varphi \in \mathcal{C}, \ \forall k\in \underline N.
\end{equation} 
 If there exists an $n$-dimensional vector $\xi_k \gg 0$ such that \eqref{cond7} holds, then for any chatter bound $N_0\geq 0$  the zero solution of the  system \eqref{DSwFDS} is GES  over the set   $ \Sigma_{\tau_a,N_0}$ of switching signals with ADT $\tau_a$ satisfying \eqref{tau}. 
\end{theorem}
\textit{Proof}. 
Let $x(t)=x(t,\sigma,\varphi)$ be the solution of \eqref{DSwFDS}, with $f_{\sigma(t)}(t,x,x_t) =f_{\sigma(t)}(t,x)+g_{\sigma(t)}(t,x_t)$,  corresponding to the initial $\varphi\in \mathcal{C}$ and the switching law $\sigma\in \Sigma_{N_0,\tau}$ having discontinuities $0<\tau_1<\ldots<\tau_p<\ldots$. Then, it follows from \eqref{metzA},\eqref{A2*} that, for any $t\in [\tau_p,\tau_{p+1})$ and  every $i\in\underline{n}, k\in \underline N$,
\begin{align}\label{mean1}
	\dot x_i(t)\ \text{sign}\;x_i(t) &=  f_{l_p,i}(t,x(t),x_t)\ \text{sign}\;x_i(t) =  f_{l_p,i}(t,x(t))\ \text{sign}\;x_i(t)\ + g_{l_p,i}(t,x_t)\ \text{sign}\;x_i(t) \notag\\
	&\leq \sum_{j=1}^n  a_{l_p,ij}(t)\;|x_j(t)| + |g_{l_p,i}(t,x_t)| ,
\end{align}
provided that the subsystem $\dot x=f_{lp}(t,x,x_t)= f_{l_p}(t,x)+ g_{l_p}(t,x_t),\  l_p \in \underline N$ is active on $[\tau_p,\tau_{p+1})$.Therefore,  it follows from \eqref{Lk}, \eqref{gk}  and \eqref{mean1} that the estimate for $D^+|x_i(t)|$ holds for all $t\in [\tau_p,\tau_{p+1}), p=0,1,...$ (with $\tau_0=0)$,
\begin{align}\label{Dini}
	D^+|x_i(t)| \leq \dot x_i(t)\;\text{sign}\;x_i(t) \leq  \sum_{j=1}^n  a_{l_p,ij}(t)\;|x_j(t)| + \sum_{j=1}^nV(\eta_{l_p,ij}	(t,\cdot))\max_{\theta\in [-h,0]}|x_j(t+\theta)|,
\end{align}
where, as defined above, for each $k\in \underline  N, \eta_{k}(t,\theta)\in NBV([-h,0],\R^{n\times n}) $ is the Riesz representation of the linear operator $L_k(t,\cdot)$. Thus, the key estimates \eqref{concl1} and \eqref{D++} in the proof of Theorem \ref{main} hold in this case. Further, the remainder of the proof is proceeded similary to that of Theorem \ref{main}, completing the proof.\hspace{13.5cm} $ \square$ 

\begin{remark}\label{all-FDE}
	{\rm
In view of the above proofs, it is important to emphasize that, if  continuous functions $A_k(\cdot), L_k(\cdot,\cdot), k\in\underline N$ are {\it given}, then Theorem \ref{main},Theorem \ref{main1} and Theorem \ref{main2}  hold true actually for \textit{any} switched  nonlinear time-varying system of the form \eqref{DSwFDS}, provided that subsystems' functions $f_k(t,x,\varphi), k\in \underline N$ satisfy the assumption (A2),(A3) or (A2*), (A3), respectively. 
}	\end{remark}

\section{Application to linear FDE systems}
In this section we will derive from the results of Section 2 some more verifiable criteria of GES for certain classes of linear FDE systems. It will be shown that, even in these particular cases, our results cover and improve a number of  results, previously known in the literature. 

First, note that in case of systems with no switchings (i.e. $N=1$) and  $f_1(t,x,x_t)= f(t,x)+g(t,x_t)$, Theorem \ref{main}  implies obviously the main result in \cite{ngoc_tinh} (Theorem 3.1) which affirms that the zero solution of the nonlinear FDE
\begin{equation}\label{tinh}
\dot x(t)= f(t,x)+g(t,x_t), \ t\geq 0,	
	\end{equation}
(where $f(t,0)=0,\ g(t,0)=0, \forall t\geq 0$) is GES if there exist continuous matrix functions $A(t)\in C(\R_+,\R^{n\times n})$ and $ \eta(t,\cdot) \in NBV([-h,0],\R^{n\times n})$ satisfying
$$  \mathcal{M}\big(J_x(f(t,x))\big)\leq A(t),\ \  |g(t,\varphi)|\leq \bigg| \int_{-h}^0 d_{\theta}[\eta(t,\theta)]\varphi(\theta)\bigg|, \ \forall t\geq 0, \ \forall \varphi \in \mathcal{C},
$$
and, moreover, there exist a vector $\xi \gg 0$ and a number $\alpha >0$ such that
$$
\big(A(t)+ e^{\alpha h} V(\eta(t,\cdot))\big)\xi \ll -\alpha \xi, \ \forall t\geq 0. 
$$ 
It is  worth noticing that the analysis in \cite{ngoc_tinh} can not be applied to deal with the problem considered in this paper since the assumption that \eqref{cond7} holds for each $k\in \underline N$  does not necessarily imply that such a {\it common} vector $\xi$ exists, as it will be  shown by the example given below. 

Further, observe that Theorems \ref{main} \ref{main1} still hold valid if one can find Meztler matrices $\widehat A _k := (\widehat a_{k,ij})\in \R^{n\times n}$ and matrix functions of bounded variations $ \widehat \eta_k(\cdot) := (\widehat\eta_{k,ij}(\cdot)) \in NBV([-h,0]),\R^{n\times n}), k\in \underline N$ with the nonnegative matrix of variation $\widehat V_k:=(\widehat v_{k,ij})= \big(Var([-h,0],\ \widehat \eta_{k,ij})\big)$ such that 
\begin{equation}\label{AB}
	\mathcal{M}(J_x(f_k(t,x,\varphi))) \leq A_k(t)\leq \widehat A_k,\ \ V(\eta_k(t,\cdot)) \leq \widehat V_k, \ \forall t \geq 0, \ \forall k\in \underline N,  
\end{equation}
and there exist vectors $\xi_k\in \R^n, \xi_k \gg 0, k\in \underline N$  such that
\begin{equation}\label{inv}
	\big(\widehat A_k+ \widehat V_k\big)\xi_k \ll 0, \   \forall k \in \underline N. 		
\end{equation}
Indeed, if \eqref{inv} holds then, by continuity, there exists a small enough $\alpha >0$ such that 
\begin{equation}\label{inv-alpha}
	(\widehat A_k+e^{\alpha h}\widehat V_k)\xi_k\ll -\alpha \xi_k,\  \forall k\in \underline N, \end{equation}
which in turn implies immediately, due to \eqref{AB}, that this set of vectors $\xi_k, \ k\in \underline N$ also satisfies \eqref{cond7}. Moreover, in this case, a  ADT lower bound $\tau_*$ in \eqref{tau} that guarantees GES of the zero solution of \eqref{DSwFDS} can be calculated explicitly for each set of vectors $\xi_k, k\in \underline N$ satisfying \eqref{inv}. Indeed, the condition \eqref{inv-alpha} can be rewritten equivalently in the form
\begin{equation}\label{cond7-i}
	g_{k,i}(\alpha) := \sum_{j=1}^n \widehat a_{k,ij}\xi_{k,j} + e^{\alpha h}\sum_{j=1}^n\widehat v_{k,ij}\xi_{k,j} + \alpha\xi_{k,i} <  0,\  i\in\underline n , \ k\in \underline N.
\end{equation}
Since, for each $k\in \underline N, i\in \underline n$, $\ g_{k,i}$ is continuous in $\alpha,\  g_{k,i}(0)<0$  and $g_{k,i}(\alpha)$ is increasing monotonically to $+\infty$ as $\alpha \rightarrow +\infty$ it implies easily that the equation $g_{k,i}(\alpha)=0$ has a unique solution $\alpha_{i,k} >0$. Setting
\begin{equation}\label{alphamax}
\alpha_{max}= \min_{k\in \underline N, \ i\in \underline n}\{\alpha_{k,i}:\  g_{k,i}(\alpha_{k,i})=0 \}
\end{equation}
it implies that  \eqref{inv-alpha} holds for all $\alpha \in [0,\alpha_{max})$ and is violated for any $\alpha \geq \alpha_{max}$.	 Hence,  $\tau_* = \frac{\ln \gamma }{\alpha_{max}}$ is the 'smallest' lower bound for ADT guaranteeing GES of the switched system \eqref{DSwFDS}. Thus, we obtain the following more applicable consequence of Theorem \ref{main}.
\begin{corollary}\label{TIS} 	Consider the switched nonlinear  FDE system \eqref{DSwFDS} which satisfies (A1,A2,A3). Assume that there exist vectors $\xi_k\in \R^n, \xi_k \gg 0, k\in \underline N$ satisfying \eqref{inv} where $\widehat A_k\in \R^{n\times n},\  \widehat \eta_k(\cdot)\in NBV([-h,0],\R^{n\times n}) $ are, respectively,  any constant matrix  and $n\times n$-matrix function of bounded variation, satisfying \eqref{AB}. Then, for each fixed chatter bound $N_0\geq 0$,   the zero solution of \eqref{DSwFDS} is GES over the set   $ \Sigma_{\tau_a,N_0}$ of switching signals with ADT $\tau_a > \tau_* = \frac{\gamma}{\alpha_{max}}$ where $\alpha_{max}$ is defined by \eqref{alphamax} and $\gamma$ is defined by \eqref{mu}.
\end{corollary} 
Now, assume that the functions $f_k(t,x,\varphi)$ are linear with respect to $x$ and $\varphi$, i.e. 
\begin{equation}\label{lin}
	f_k(t,x,\varphi) = A_k(t)x + L_k(t,\varphi),\ t\geq 0, \ k\in \underline N, 
\end{equation}
 where $A_k(t)$ and $L_k(t,\varphi)$ satisfy the continuity assumptions as in Theorem \ref{main}. Then, since $J_x(f_k(t,x,\varphi))=A(t)x, \ f_k(t,0,\varphi)= L_k(t,\varphi)$, from Theorem  \ref{main} we get  straightforwardly the following result which can be considered as an extension of the  main results proved for non-switching linear systems in \cite{ngoc_tinh} (Corollary 3.3) and \cite{ngoc_tinh_huy} (Theorem 3.2) to the class of switched linear FDE systems, with average dwell time switchings.
\begin{corollary}\label{linear}
	Consider the time-varying switched linear  FDE system of the form
	\begin{equation}\label{sw_LFDE}
		\dot x(t)= A_{\sigma(t)}(t)x(t)+ \int_{-h}^0d_{\theta}[\eta_{\sigma(t)}(t,\theta)]x(t+\theta),\ t\geq 0,\  \sigma \in \Sigma_+.
	\end{equation}
	Assume that there exist vectors $\xi_k\in \R^n, \xi_k \gg 0, k\in \underline N$ and a positive number $\alpha >  0$ such that
	\begin{equation} 
		\label{cond7-l}
		\big(\mathcal{M}\big(A_k(t)\big)+ e^{\alpha h} V\big(\eta_k(t,\cdot)\big)\big)\xi_k \ll -\alpha \xi_k, \ \forall t\geq 0,\  \forall k \in \underline N. 
	\end{equation}
	Then the  switched linear system \eqref{sw_LFDE} is GES over the set of switching signals  $\Sigma_{\tau_a}$ with ADT $\tau_a > \tau_* = \frac{\gamma}{\alpha}$, where $\gamma$ is defined by \eqref{mu}.
\end{corollary}
It is important to note that the matrices of variations $V(\eta_k(t,\cdot))$ in \eqref{cond7} and \eqref{cond7-l} can be explicitly calculated in some  particular cases of interest. For instance, assume that the matrix functions of bounded variations $\eta_k(t,\cdot), k\in \underline N, $ in the Assumption (A3) are  given by
\begin{equation}\label{etak}
	\eta_k(t,\theta)= \sum_{i=1}^mB_k^{i}(t)\chi_{(-h^i_k(t),0]}(\theta)+\int_{-h}^{\theta}C_k(t,s)ds, \ t\geq 0, \ \theta\in [-h,0],
\end{equation}
 where, for each  $ k\in \underline N, i\in \underline m, B^i_k(t), C_k(t,\theta)$ are continuous $(n\times n)$-matrix functions,  $ h\geq h^{m}_k(t)>\ldots> h^2_k(t)>h^1_k(t)>0, \forall t\geq 0 $, with $h^i_k(\cdot)$ being given continuous functions of $t$ and $ \chi_M$ is, by definition, the characteristic function of the set $M\subset \R$.   Note that, in this case, the associated linear operators $L_k$ defined by \eqref{Lk} are given by
 \begin{equation}\label{Lphi}
 	L_k(t,\varphi) = \sum_{i=1}^m B^i_k(t)\varphi (-h_k^i(t))+ \int_{-h}^{0}C_k(t,s)\varphi (s) ds, \ t\geq 0, \ \varphi \in \mathcal{C}.
 \end{equation}Then, by a simple calculation, we get 
 \begin{equation}\label{Vk}
 	V(\eta_k(t,\cdot))\ \leq \ \sum_{i=1}^m|B_k^{i}(t)|+\int_{-h}^{0}|C_k(t,s)|ds,
 \end{equation}	
where the equality holds, provided that $B^i_k(t)\geq 0, C_k(t,\theta) \geq 0, \forall t\geq 0, \forall\theta \in [-h,0]$. 
Moreover, in this case the linear system \eqref{sw_LFDE} is reduced to the switched linear system with multiple discrete delays and distributed delay of the form
\begin{equation}\label{LDE}
	\dot x(t)= A_{\sigma(t)}(t)x(t)+ \sum^m_{i=1}B^i_{\sigma(t)}(t)x(t-h_{\sigma(t)}^i(t))	+ \int_{-h}^{0}C_{\sigma(t)}(t,\theta)x(t+\theta)d\theta, \ t\geq 0. 
	\end{equation}
Therefore, as an immediate consequence of Corollary \ref{linear}, we have
\begin{corollary}\label{discretedelay} 
	Consider the time-delay switched linear  system \eqref{LDE}. Assume that there exist vectors $\xi_k\in \R^n, \xi_k \gg 0, k\in \underline N$ and a positive number $\alpha >  0$ such that
	\begin{equation} 
		\label{cond9}
		\mathcal{M}\big(A_k(t)\big)\xi_k + e^{\alpha h}\bigg(\sum_{i=1}^m|B_k^{i}(t)|+\int_{-h}^{0}|C_k(t,s)|ds\bigg) \xi_k \ll -\alpha \xi_k, \ \forall t\geq 0,\  \forall k \in \underline N. 
	\end{equation}
	Then the switched linear system \eqref{LDE} is GES over the set of switching signals  $\Sigma_{\tau_a}$ with ADT $\tau_a > \tau_* = \frac{\gamma}{\alpha}$, where $\gamma$ is defined by \eqref{mu}.
	\end{corollary}

Applying Corollary \ref{TIS}-Corollary \ref{discretedelay} to the system \eqref{LDE} we get  the following

\begin{corollary}\label{Qi}
	Consider the time-delay switched linear  system \eqref{LDE}. Assume that there exist constant matrices $ \widehat A_k, \widehat B^i_k \in \R^{n\times n}, i\in \underline m, k\in \underline N$ and continuous $(n\times n)$-matrix functions $ \widehat C_k(\theta), k\in \underline N,\ \theta\in[-h,0]$ such that, for each $k\in \underline N, t\geq 0$ we have $  A_k(t)\leq \widehat A_k, |B_k^i(t)|\leq \widehat B_k^i, \forall i\in \underline m$ and $|C_k(t,\theta)|\leq \widehat C_k(\theta), \theta \in [-h,0]$. Assume moreover that there exist positive vectors $\xi_k \gg 0, k\in \underline N,$ such that
	\begin{equation}\label{hat1}
	\big(\mathcal{M}(\widehat A_k) + \widehat V_k\big) \xi_k \ll 0, \ \forall k \in \underline N, 
	\end{equation}
	where $ \widehat V_k:= \sum_{i=1}^m \widehat B_k^{i}+\int_{-h}^{0}\widehat C_k(s)ds$. Then the time-delay switched linear system \eqref{LDE} is GES over the set of switching signals  $\Sigma_{\tau_a}$ with ADT $\tau_a >\tau_*= \frac{\gamma}{\alpha_{max}}$ where $\alpha_{max}$ is defined by \eqref{cond7-i}, \eqref{alphamax} and $\gamma$ is defined by \eqref{mu}.
\end{corollary} 

If vectors $\xi_k, k\in \underline N$ in the above Corollaries \ref{TIS} -  \ref{Qi} can be chosen the same  $\xi_k =\xi, k\in \underline N$, then by Theorem \ref{main1} we obtain, correspondingly, the sufficient conditions for GES of the switched nonlinear system  \eqref{DSwFDS} and the switched  linear systems \eqref{sw_LFDE}, \eqref{LDE} over the set of switching signals $\Sigma_+$. For instance, by Theorem \ref{main1} and Corollary \ref{TIS}, Corollary \ref{discretedelay} we get
\begin{corollary}\label{identxi}
The time-delay switched linear system \eqref{LDE} is GES over the set of switching signals  $\Sigma_{+}$ if there exist a positive vector $\xi\in \R^n, \xi \gg 0,$ and $\alpha >0 $ such that
\begin{equation} 
	\label{cond10}
\bigg(	\mathcal{M}\big(A_k(t)\big)\xi + e^{\alpha h}\sum_{i=1}^m|B_k^{i}(t)|+ e^{\alpha h}\int_{-h}^{0}|C_k(t,s)|ds\bigg) \xi \ll -\alpha \xi, \ \forall t\geq 0,\  \forall k \in \underline N. 
\end{equation}
In particular, if the system \eqref{LDE} is time-invariant, i.e. $A_k(t)\equiv A_k, B_k(t)\equiv B_k, C_k(t,\cdot)\equiv C_k(\cdot), \ \forall t \geq 0,\ \forall k\in \underline N$ then the above condition can be replaced by
\begin{equation}\label{hat2}
\big(	\mathcal{M}\big(A_k\big) +  V_k \big) \xi \ll 0, \ \forall k\in \underline N,
\end{equation}
where $V_k:= \sum_{i=1}^m|B_k^{i}|+ \int_{-h}^{0}|C_k(s)|ds$. 
	\end{corollary}

Particularly, when $m=1$ and $C_k(t,s)\equiv C_k(t), \forall s\in [-h,0]$, Corollary \ref{Qi} is reduced to Theorem 1 in \cite{QiSun}, see also Corollary 3.1 in \cite{Tian_Sun} and Corollary 3.2 \cite{Li2018}. Similarly, if all the assumptions of Corollary \ref{Qi} hold and, moreover, $\xi_k=\xi, \ \forall k\in \underline N$ then from Corollary \ref{identxi} we get back to the main results in \cite{Sun2012, Li2017}.

\begin{remark}\label{A_0} {\rm 
	It is well-known that the condition  \eqref{hat1}) is equivalent (see, e.g. \cite{HornJohn} and \cite{Son_Hin96}) to that the Metzler matrices $\widehat B_k:=\mathcal{M}(\widehat A_k)   + \widehat V_k,\  k\in \underline N $ are all \textit{Hurwitz stable } (i.e. all zeros $\lambda$  of the characteristic polynomials $P_k(\lambda):= \det (\lambda I-\widehat B_k), k\in \underline N	$ have negative real parts: ${\rm Re }\lambda <0$). Moreover, in order to verify the existence of a common vector $\xi \gg 0$ satisfying $ (\mathcal{M}(A_k) +V_k)\xi \ll 0, \ \forall k\in \underline N$ as required in Corollary \ref{identxi} we can apply the procedure described in \cite{Knorn_Mason2009} to Metzler matrices $B_k:=\mathcal{M}(A_k) +V_k, k\in \underline N$ .  }	 \end{remark}

The approach developed in this paper is applicable to study exponential stability of \textit{switched positive FDE systems}. Below, we demonstrate this applicability  for the  case of time-invariant linear switched FDE systems.

It is immediate from Definition \ref{positive} that the switched system \eqref{DSwFDS} is positive if and only if all its subsystems  \eqref{DSwFDSk} are positive. In particular, it follows that the time-invariant switched linear FDE system 
\begin{equation}\label{invLFDE}
	\dot x(t)= A_{\sigma(t)}x(t)+ \int_{-h}^0d_{\theta}[\eta_{\sigma(t)}(\theta)]x(t+\theta),\ t\geq 0, \ \sigma \in \Sigma,
\end{equation}
is positive if and only if  all $A_k, k\in \underline N$ are Metzler matrices and all the functions $\eta_k(\cdot) \in NBV([-h,0],\R),$ $ k\in \underline N$ are non-decreasing on $[-h,0]$, i.e. $ \eta_k(\theta_1) \leq \eta_k(\theta_2) $ whenever $-h\leq \theta_1<\theta_2\leq 0$, see, e.g. \cite{ngocFDE}. The following result gives a verifiable criterion for the switched positive linear  system  \eqref{invLFDE} to be  GES under arbitrary switchings with ADT.
\begin{corollary}\label{stabpositive}
Consider the time-invariant switched positive linear FDE system \eqref{invLFDE}. Assume that there exist vectors $\xi_k\in \R^n, \xi_k \gg 0, k\in \underline N$ such that
\begin{equation}\label{GESinv}
	(A_k+\eta_k(0))\xi_k \ll 0, \ \forall k\in \underline N.
\end{equation}
	
\noindent Then the system \eqref{invLFDE} is GES over the set of switching signals $\Sigma_{\tau_a}$ satisfying $\tau_a >\tau_*= \frac{\ln\gamma}{\alpha_{max}}$  where $\gamma$ is defined by \eqref{mu} and $\alpha_{max}:= \min_{k\in \underline N,i\in \underline n}\alpha_{k,i}$,  with  $\alpha_{k,i}$ being defined as the solutions of the equations 
\begin{equation} \sum_{j=1}^n\big(a_{k,ij} + e^{\alpha h} \eta_{k,ij}(0)\big)\xi_{k,j}+ \alpha\xi_{k,i}=0, \ k\in \underline N, \  i\in \underline n.
\end{equation} 
If, moreover, there exists $\xi \gg 0$ such that \eqref{GESinv} is satisfied for $\xi_k=\xi, \forall k\in \underline N$ then the system \eqref{invLFDE} is GES over the set of switching signals $\Sigma_+$. 
Conversely, if the switched positive linear system \eqref{invLFDE} is is GES over a set of switching signals $\Sigma_{\tau_a}$ with $ \tau_a >0$ then there exist vectors $\xi_k\in \R, \xi_k\gg 0, k\in \underline N$ such that  \eqref{GESinv} is satisfied. 
	\end{corollary} 

{\it Proof. } 
Since the matrix functions $\eta_k(\theta)$ are non-decreasing for $\theta \in [-h,0]$ and $\eta_k(-h)=0$, by the definition of $NBV([-h,0],\R^{n\times n})$, we have obviously $V(\eta_k)=\eta_k(0)$, for each $k\in \underline N$. Therefore, the first part of the Corollary \ref{stabpositive} is immediate from Corollary \ref{linear} and Theorem \ref{main1}. Conversely, if \eqref{invLFDE} is positive and GES over  $\Sigma_{\tau_a}, \tau_a >0$ then the linear  subsystem $(A_k,\eta_k)$ is positive and GES, for each $k\in \underline N$. The latter is equivalent to that the Metzler matrices $A_k+\eta_k(0), k\in \underline N$ are Hurwitz stable, by Theorem 4.1 of \cite{ngocFDE}, which in turn is equivalent to the existence of vectors $\xi_k \gg 0, k\in \underline N$ satisfying \eqref{GESinv} (see, e.g. \cite{Son_Hin96}). This completes the proof.\hskip5.0cm $\square$

Corollary \ref{stabpositive}, in such a general formulation,  is  novel in the literature on stability of switched positive systems. Specifically, we can apply Corollary \ref{stabpositive} to get a similar criterion of GES for time-invariant switched positive systems with mixed delays of the form \eqref{LDE} where, for each $k\in \underline N$, $A_k(t)\equiv A_k$ are Metzler matrices, $B^i_k(t)\equiv B^i_k\geq 0$ and $C_k(t,s)\equiv C_k(s)\geq 0, \forall s\in [-h,0]$. It is notable that, even in this particular case, criteria of global exponential stability for switched positive linear systems using ADT switching signals are still lack so far in the literature. 

Finally, it is worth noticing that Corollary \ref{stabpositive}, in the case of non-delay systems (i,e, $ \eta_k=0, \forall k$), implies Corollary 3.6. and Theorem 3.8 in \cite{Dong2015}. Also,  the main result of \cite{Liu_Dang}  (Theorem 2) is obviously followed from Corollary \ref{stabpositive},  by assuming $\xi_k=\xi, \ k\in \underline N$,  and choosing $\eta_k(\cdot)$ to be as \eqref{etak}, with $B^i_k(t)\equiv B^i_k\geq 0, \ i\in \underline m, \ C_k\equiv 0, \forall k\in \underline N$,  which yields readily  $ \eta_k(0)= \sum_{i=1}^m B^i_k$.


\section{Application to nonlinear time-delay systems with sector nonlinearity}

In this section we will show that the approach developed in Section 2 can be applied to get sufficient conditions of absolute exponential stability for a class of switched time-delay systems with sector nonlinearity. 

Consider the time-varying nonlinear switched system of the form
\begin{equation}\label{sector}
\dot x(t)=P_{\sigma(t)}(t)\psi (x(t))+ B_{\sigma(t)}(t)\psi (x(t-h)),\ t\geq 0,\ \sigma\in \Sigma_+,
\end{equation}
where the set of switching signals $\Sigma_+$ is defined as in the Section 2. The corresponding family of constituent subsystems is given by
\begin{equation}\label{sector_sub}
	\dot x(t)=P_k(t)\psi (x(t))+ B_k(t)\psi (x(t-h)),\ t\geq 0,\ k=1,\ldots, N.
\end{equation}
Here $P_k(\cdot), B_k(\cdot)$ are continuous $(n\times n)$-matrix functions, and the nonlinearity $\psi:\R^n\rightarrow \R^n $ is assumed to be continuous, diagonal
$$
\psi (x)=(\psi_1(x_1) \ \psi_2(x_2)\ \ldots \ \psi_n(x_n))^{\top}
$$
and satisfy
\begin{equation}\label{sector_positive}
0<	x_i\psi_i(x_i) \ \text{for}\ x_i\not=0,\ i\in \underline n,
\end{equation}
or, equivalently, $\text{sign}\;x_i\ \psi(x_i) > 0$ for $x_i\not=0$. Such a nonlinearity is said to be \textit{admissible}.  It follows immediately from the continuity that $\psi_i(0)=0, \ \forall i\in \underline n$ for any admissible $\psi $ and therefore the system \eqref{sector} admits the zero solution $x(t)\equiv 0, \ t\geq 0$. Moreover, without speaking additionally, the function $\psi $ is  assumed to satisfy certain
 Lipschitz condition that ensures  global existence and  uniqueness of the solution of system \eqref{sector}, for any initial condition $\varphi \in \mathcal{C}$ and any switching $\sigma \in \Sigma_+$.

 \begin{definition} System \eqref{sector} is said to be absolutely exponentially stable (shortly AES) over  the set of switching signals $\Sigma_+$ if its zero solution is GES for any switching signal $\sigma \in \Sigma_+$ and any admissible nonlinearity $\psi(x)$.
 	\end{definition}
 Dynamics with sector nonlinearities  are widely used in modeling automatic control systems and neural networks (see, e.g. \cite{KaszBhaya}). In case of non-switched systems the conditions of absolutely asymptotic stability  were investigated in \cite{persid69},\cite{KaszBhaya}  by using the Lyapunov direct method. More recently, similar problems has been studied for  switched time-delay systems of the form \eqref{sector} in a number of  works, but only for the  time-invariant case, i.e. when $P_k(t)\equiv P_k, B_k(t)\equiv B_k, \forall t$), mostly by using  Lyapunov - Krasovskii functionals (see, e.g. \cite{alex_platon2008, sun_wang2013, alex_mason2014,  alex_zab2017},  and also \cite{zhang_han2014, zhang_zhao2016} for conditions of absolute exponential stability).  Below, we will use Theorem \ref{main2} to derive a new criterion of absolutely exponential stability for the time-varying switched system \eqref{sector}  with ADT switchings.
 
\begin{theorem}\label{ab1}
Consider the  switched system \eqref{sector} with an  admissible nonlinearity $\psi$ satisfying the sector constraint of the form
\begin{equation}\label{sectork}
0 <  x_i \psi_i(x_i)\ \leq \ \beta_i x_i^2, \ \text{\rm for} \ x_i\not=0, \ \  i=1,2,\ldots, n, \ 
\end{equation}
where $ \beta_i  > 0 $ are given positive numbers. Assume that there exists an $n$-dimensional vector $\xi_k \gg 0$ and a real number $\alpha >0$ such that 
\begin{equation}\label{cond12}
(\mathcal{M}(P_k(t))	+ e^{\alpha h}|B_k(t)|)D_{\beta}\xi_k \ll -\alpha \xi_k,\ \forall t \geq 0, \ \forall k\in \underline N,
	\end{equation}
 where $D_{\beta}$ is a diagonal matrix defined as 
 $$
 D_{\beta}= \text{\rm diag}(\beta_1, \beta_2, \ldots, \beta_n).
 $$Then for any chatter bound $N_0\geq 0$ the zero solution of the  system \eqref{sector} is AES over the set   $ \Sigma_{\tau_a,N_0}$ of switching signals with ADT $\tau_a > \tau_*= \frac{\ln \gamma}{\alpha}$  where $\gamma $ is defined by \eqref{mu}. Moreover, if there exists a vector $\xi \gg 0$ such that \eqref{cond12} holds for all $ \ \xi_k=\xi,\ k\in \underline N$ then the zero solution of \eqref{sector} is AES over the set of switching signals $\Sigma_+$.
\end{theorem} 
\textit{Proof}. First, it is straightforwardly verified that for any admissible nonlinearity satisfying \eqref{sectork} we have, for each $i\in \underline n, $
\begin{equation}\label{psi}
		\psi_i(x_i)\; \text{sign} x_i = |\psi_i(x_i)|\leq \beta_i |x_i|,\ \text{and}\ \ \psi_j(x_j)\;\text{sign}\; x_i \leq |\psi_j(x_j)|, \ \ \forall x_i\in \R,\  \forall j\in \underline n, j\not=i,
\end{equation}
(where the case $x_i\not=0$ is implied from \eqref{sectork} while the case $x_i=0$ is immediate because $\psi_i(0)=0$). For each $k\in \underline N$, define the continuous functions $f_k(t,x)\ =\! P_k(t)\psi (x)$ and $  g_k(t,\varphi)\!=\! B_k(t)\psi (\varphi(-h))$ with $ t\geq 0,\  x\in \R^n,\ \varphi \in \mathcal{C}. $  Then, we deduce readily, by using \eqref{psi}, 
that, for each  $k\in \underline N$ and $i\in \underline n$, 
\begin{align*}
f_{k,i}(t,x) \ \text{sign}\ x_i &=  p_{k,ii}(t)\psi _i(x_i)\ \text{sign}\ x_i  +  \sum_{ j=1, j\not=i}^n p_{k,ij}(t)\psi_j(x_j)  \text{sign}\ x_i \\
&\leq p_{k,ii}(t) |\psi _i(x_i) | + \sum_{j=1,j\not=i }^n |p_{k,ij}(t)|\ |\psi_j(x_j)|\\
& \leq p_{k,ii}(t)\beta_i |x_i | + \sum_{j=1,j\not=i }^n |p_{k,ij}(t)|\beta_j |x_j| = \sum_{i=1}^n a_{k,ij}(t)|x_j|, \ \forall t\geq 0, \ \forall x\in \R^n, 
\end{align*}
where, by definition, $A_k(t)= (a_{k,ij}(t)) := \mathcal{M} (P_k(t)) D_{\beta}$, being  clearly a Metzler matrix. Thus, $f_k(t,x), k\in \underline N,$ satisfy the assumption (A2*).
On the other hand, define for each $t\geq 0 $ and $k\in \underline N$ the  linear operators $L_k(t,\cdot):  \mathcal{C} \rightarrow \R^n$ by setting 
\begin{equation}\label{A_Lk}
	L_k (t,\varphi)= |B_k(t)|D_{\beta}\varphi(-h),\  \varphi \in \mathcal{C}.
\end{equation} 
Then,  $L_k(t,\cdot) $ are positive operators which have the corresponding Riesz representation of the form \eqref{Lk} with  $\eta_k(t,\cdot) \in NBV([-h,0],\R_+^{n\times n})$ defined by
\begin{equation}\label{etak2}
	\eta_k(t,\theta)=\begin{cases}  0 \ &\text{if}\ \ \theta=-h,  \ \forall  t\geq 0 \\  |B_k(t)|D_{\beta} \ &\text{if}\ \  \theta \in (-h,0], \ \forall t \geq 0  \end{cases}
\end{equation} 
	Then, for each addmissble nonlinearity $\psi$ satisfying \eqref{sectork}, we get, by using \eqref{psi} and the positivity of $L_k(t,\cdot)$,
\begin{align*}
|g_k(t,\varphi)|&\leq |B_k(t)||\psi(\varphi(-h))|\stackrel{\eqref{psi}}\leq |B_k(t)|D_{\beta} |\varphi(-h)|\\
&= L_k(t,|\varphi|) = |L_k(t,\varphi)|,\ \forall t \geq 0, \ \forall \varphi \in \mathcal{C}, \ \forall k\in \underline N.
\end{align*}
Thus, \eqref{gk}  is also satisfied. Moreover, we have obviously $V(\eta_k(t,\cdot)):= \big(Var([-h,0], \eta_{k,ij}(t,\cdot))\big)=|B_k(t)|D_{\beta}$. Therefore, the assertions of the theorem are followed straightforwardly from Theorem \ref{main2} and Theorem \ref{main1}, taking into account Remark \ref{all-FDE}, completing the proof. \hspace{15.8cm}$\square$

\begin{remark}\label{7}
	{\rm  It is important to emphasize that conditions for absolute stability of switched systems with sector nonlinearities have been so far obtained only in the time-invariant case (see, e.g. \cite{persid69, persid2005, alex_platon2008, sun_wang2013,   alex_mason2014, zhang_zhao2016} ) To the best of our knowledge, Theorem \ref{ab1} seems to be the first result in the existing literature which gives the condition of absolute exponential stability for switched  time-varying systems. The other distinction between our above result and the results in the mentioned previous works is that our proofs are based on the comparison principle while all the latters have made use of the Lyapunov - Krasovskii functionals as a main tool. Moreover, our approach can be applied  to more general systems, in the unified way. In particular, the case of systems with multiple discrete delays and distributed delays, rather than that with a single delay of the form \eqref{sector} can also be treated similarly, by the above approach. }
\end{remark}
	
\begin{remark}\label{8}{\rm 
	Even in the time-invariant situation our above result extends some previously known results and is reduced to less conservative criteria. To see that, assume $P_k(t)\equiv P_k,\  B_k(t)\equiv B_k, \forall t\geq 0$ then it follows from Theorem \ref{ab1} (by letting  $\alpha \downarrow 0$ and $\zeta_k := D_{\beta}^{-1}\xi_k $ in \eqref{cond12}) that the existence of a vector $\zeta \gg 0$ satisfying  $N$ systems of linear inequalities
	\begin{equation}\label{sun}
	(\mathcal{M}(P_k) + |B_k|)\zeta \ll 0,\  k\in \underline N	
	\end{equation}
    is sufficient to ensure that the switched nonlinear system 
	$$
	\dot x(t)= P_{\sigma(t)} \psi(x(t)) + B_{\sigma(t)}\psi(x(t-h)), \ t\geq 0, \ \sigma\in \Sigma_+	
	$$  is absolutely exponentially stable, provided that the admissible sector nonlinearity $\psi$ satisfies the constraint \eqref{sectork}. The condition \eqref{sun}, in some sense, is less restrictive and easier to be checked than those given in previous works which are based on building up a common Lyapunov- Krasovskii for the system. For instance, the condition  in \cite{sun_wang2013} (Theorem 2) requires $\zeta \gg 0$ to satisfy  $N$  {\it dual systems} of linear inequalities
	 \begin{equation}\label{sw}
		(\mathcal{M}(P_k) + \widetilde B_k)^{\top}\zeta \ll 0, \ k\in \underline N, \text{with}\  \widetilde B_k= \big(\max_{k\in\underline N}|b_{k,ij}|\big)
	\end{equation}
 while in \cite{alex_mason2014} (Theorem 3.1), such a $\zeta$ has to satisfy  $N^2$  dual systems 
\begin{equation}\label{alexmas}
	(\mathcal{M}(P_k) + |B_s|)^{\top}\zeta \ll 0,\  k,s \in \underline N. 
\end{equation} 
 These conditions are obviously more restrictive than \eqref{sun} as shown by an example in the next section. This is even more evident in the case of non-delay positive systems of the form \eqref{sector} (i.e. when $\mathcal{M}(P_k) =P_k, B_k=0$): while, by our result, the existence of strictly positive solution $\zeta$ of $P_k\zeta \ll 0, k\in \underline N$ implies absolute exponential  stability, the well-known result of \cite{alex_platon2008} (Theorem 3) requires such a condition to hold for the both original and dual inequalities: there exist $\zeta\gg, \nu \gg 0$ such that
		$$ 	P_k\zeta \ll 0 \ \text{and}\  P_k^{\top}\nu \ll 0, \  \forall k\in \underline N.	$$   
		It is worth mentioning moreover that, differently from the previous works, the approach of this paper, based on ADT switchings, allows us to deal with the situation when such a common positive vector $\zeta$ does not exist, proving not only absolute exponential stability of the system but also giving an estimate on the convergence rate of the solution to the equilibrium.  However, comparing with   \cite{alex_platon2008, sun_wang2013} and \cite{alex_mason2014, alex_zab2017}, we have to restrict the class of sector nonlinearities by the condition \eqref{sectork}.  }
	\end{remark}

\begin{remark}\label{9}{\rm  Obviously, the class of admissible nonlinearities  \eqref{sectork} includes the class of nonlinearities $\psi$ satisfying
		$$
	\delta_ix_i^2 \leq x_i\psi_i(x_i) \leq \beta_i x_i^2,\ \text{for} \ x_i\not= 0, \ \forall i\in \underline n
		$$
		(with $\beta_i\geq \delta_i > 0, i\in \underline n$ being given positive numbers) which has been considered in \cite{zhang_han2014, zhang_zhao2016}, where a similar problem is studied for time-invariant switched  systems with ADT switchings, also by using the Lyapunov- Krasovskii functional. It is easily verified that the above Theorem \ref{ab1}, when applies to this case  is reduced to the condition that is less conservative than the main results of \cite{zhang_zhao2016} (Theorem 1, Theorem 2).
}
\end{remark}

\section{Illustrative examples}
We give first an example to illustrate the use of Theorem \ref{main} in the case where $f_k(t,x,x_t)$ can not be represented in the form \eqref{tinh}. 

\begin{example}	\label{ex1}
	{\rm Consider a nonlinear switched system of the form \eqref{DSwFDS} in $\R^2 $ with $N=2, h=1$ where, for all $t\geq 0,$
		\begin{equation*}
			f_1(t, x,x_t)=
			\begin{bmatrix}
				-6x_1+\sqrt{x_1^2+x_2^2}. \sin^2 t +2\sin t.x_1(t-1)+\cos t.x_2(t-1)\\
				x_1.\sin^2(x_2(t-1))-5x_2 +x_2.\sin^2(x_1(t-1))+\cos t.x_1(t-1) +2\cos t.x_2(t-1)
			\end{bmatrix},
		\end{equation*}	
		\begin{equation*}
			f_2(t, x,x_t)=
			\begin{bmatrix}
				-5x_1+x_1.\cos^2(x_1(t-1))	+x_2.\sin^2(x_2(t-1)) +2\cos t.x_1(t-1) +\cos t.x_2(t-1)\\
				\sqrt{x_1^2+x_2^2}. \cos^2 t-6x_2+\sin t.x_1(t-1)+2\sin t.x_2(t-1)
			\end{bmatrix}.
		\end{equation*}			
		A simple calculation shows that
		\begin{equation*}
			\mathcal M\big(J_xf_1\big)\!=\!	\begin{bmatrix}
				-6+\frac{x_1}{\sqrt{x_1^2+x_2^2}} .\sin^2 t\! &\frac{|x_2|}{\sqrt{x_1^2+x_2^2}} .\sin^2 t \\
				\sin^2(x_2(t-1)) \!&-5+\sin^2(x_1(t-1))
			\end{bmatrix}\leq A_1(t):=
			\begin{bmatrix}
				-6+\sin^2 t&\sin^2 t \\
				1&-4
			\end{bmatrix}\leq \widehat A_1:= \begin{bmatrix}-5&1\\1&-4\end{bmatrix},
		\end{equation*}		
		\begin{equation*}
			\mathcal M\big(J_xf_2\big)\!=\!		\begin{bmatrix}
				-5\!+\cos^2(x_1(t-1))  &\sin^2(x_2(t-1))\\ \frac{|x_1|}{\sqrt{x_1^2+x_2^2}} .\cos^2 t &\!\!-6\!+\!\frac{x_2}{\sqrt{x_1^2+x_2^2}} .\cos^2 t
			\end{bmatrix}\!\leq A_2(t):=
			\begin{bmatrix}
				-4&	1 \\
				\cos^2 t&-6\!+\!\cos^2 t
			\end{bmatrix}\leq \widehat A_2:= \begin{bmatrix}-4&1\\1&-5\end{bmatrix},
		\end{equation*}
		and, in view of \eqref{Vk}, 	
		\begin{equation*}
			V(\eta_1(t,\cdot)) \leq
			\begin{bmatrix}
				2|\sin t|&|\cos t|\\
				|\cos t|&2|\cos t|
			\end{bmatrix} \leq \widehat V_1:=\begin{bmatrix}2&1\\1&2\end{bmatrix},\
			V(\eta_2(t,\cdot))\leq
			\begin{bmatrix}
				2|\cos t| &|\cos t|\\
				|\sin t|&2|\sin t|
			\end{bmatrix} \leq \widehat V_2:=\begin{bmatrix}2&1\\1&2\end{bmatrix}.
		\end{equation*}
		Then, taking  $t=0$ and $\alpha=0$, we have
		\begin{equation*}
			A_1(0)+V(\eta_1(0,\cdot))=
			\begin{bmatrix}
				-6&	0 \\
				1&-4
			\end{bmatrix}+\begin{bmatrix}
				0&1 \\
				1&2
			\end{bmatrix}=
			\begin{bmatrix}
				-6&1 \\
				2&-2
			\end{bmatrix}
		\end{equation*}
		\begin{equation*}
			A_2(0)+V(\eta_2(0,\cdot))=
			\begin{bmatrix}
				-4&	1 \\
				1&-5
			\end{bmatrix}+\begin{bmatrix}
				2&1 \\
				0&0
			\end{bmatrix}=
			\begin{bmatrix}
				-2&2 \\
				1&-5
			\end{bmatrix},
		\end{equation*}
		and therefore, there does not exist a vector $\xi = [c_1,c_2]^{\top}>0$ such that $\big(A_k(0)+ V(\eta_k(0,\cdot))\big)\xi \ll 0,$ $k=1,2,$ because otherwise, we would get $c_1<c_2$
		and $c_2 < c_1$,  a contradiction. Therefore, Theorem \ref{main1} does not apply in this case. On the other hand, it is obvious that 
		\begin{equation*}
			A_1(t)+V(\eta_1(t,\cdot))\leq \widehat A_1+\widehat V_1=
			\begin{bmatrix}
				-3&2 \\
				2&-2
			\end{bmatrix}, \ \  A_2(t)+V(\eta_2(t,\cdot))\leq \widehat A_2+\widehat V_2= 
			\begin{bmatrix}
				-2&2 \\
				2&-3
			\end{bmatrix},\ \forall t \geq 0, 
		\end{equation*}
	and it can be easily verified that  \eqref{inv} holds for  the two vectors $\xi_1 = [0.8\ 1]^{\top}$ and $\xi_2 = [1\  0.8]^{\top}$,  with $\alpha=0.1013$.  Therefore, by Corollary \ref{TIS}, we conclude  that the switched nonlinear FDE under consideration is GES over the set $\Sigma_{\tau_a}$ of switching signals with ADT $\tau_a >\frac{\ln \gamma}{\alpha}=\frac{\ln 1.25}{0.1013}\approx 2.2028.$ For instance,  if we  choose the switching signal $\sigma\in \Sigma_{\tau_a}$ with  $\tau_a=3>\tau_*=2.2028$, as shown in Fig. 1 and the initial condition given by the function $\varphi(\theta)\equiv (-1\; \cos \theta)^{\top},\;\theta\in [-1, 0]$, then the solution's trajectory of the above nonlinear switched system converges exponentially to zero, as shown in Fig. 2. The simulation has been performed with the MATLAB code {\bf dde23}.
		\begin{figure}[ht]
			\begin{center}
				\includegraphics[height=6cm,width=11cm]{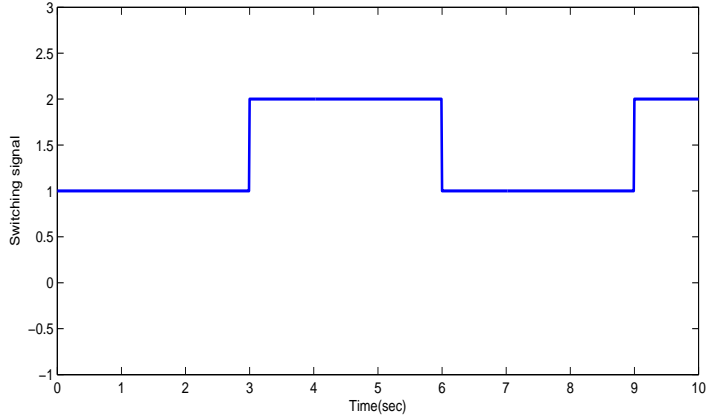}
			\end{center}\vspace{-0.5cm}
			\caption{\textit{The switching with ADT $\tau_a=3$}}
		\end{figure}
		
	\begin{figure}[ht]
			\begin{center}
				\includegraphics[height=6cm,width=11cm]{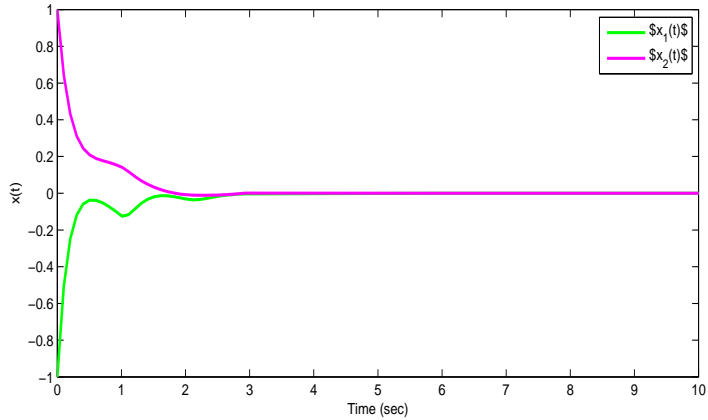}
			\end{center}\vspace{-0.5cm}
		\caption{\textit{The  solution trajectory of Example 1  under switching with ADT $\tau_a=3>2.2028$}}
		\end{figure}
	}
\end{example}
The next example illustrates the use of Theorem \ref{ab1} to affirm the absolute exponential stability of a time-delay system of the form \eqref{sector}. 
\begin{example}\label{ex2}{\rm 
		Consider a system \eqref{sector} with $n=2, N=2, h=1$ and
	\begin{equation*}
	P_1= \begin{bmatrix}-2&0\\1 &-4
		\end{bmatrix}, \  P_2= \begin{bmatrix}-2&0\\0 &-2
		\end{bmatrix}, \ B_1=\begin{bmatrix} 0&2\\0&2
\end{bmatrix},\ \ B_2= \begin{bmatrix}1&1\\1 &0
\end{bmatrix}, \ \widetilde B: = \big(\max_{k\in\underline N}|b_{k,ij}|\big)= \begin{bmatrix} 1&2\\1&2
\end{bmatrix}.
\end{equation*}	 We have 
\begin{equation*}
	P_1+B_1=\begin{bmatrix} -2&2\\1&-2
	\end{bmatrix}, \ P_2+B_2= \begin{bmatrix} -1&1\\1&-2 \end{bmatrix}, \ \ P_1+\widetilde{B}=\begin{bmatrix} -1&2\\2&-2\end{bmatrix},\  P_2+\widetilde{B}=\begin{bmatrix} -1&2\\1&0\end{bmatrix}.
\end{equation*}
It is clear that there is no $\zeta \gg 0$ satisfying \eqref{sw} and hence Theorem 2 of \cite{sun_wang2013} can not be applied in this case. Also, it can be checked that there is no $\zeta \gg 0$ such that $(P_k+B_k)^{\top}\zeta \ll 0, \ k=1,2$ and thus \eqref{alexmas} fails to hold. Therefore, Theorem 3.1 of \cite{alex_mason2014} cannot be used to construct a Lyapunov-Krasovskii functional for this system. On the other hand, it is easily verified that vector $\zeta= (7 \ 4)^{\top}$ satisfies \eqref{sun}, so that we can apply Theorem \ref{ab1} to affirm that the switched time-delay system under consideration is AES, for any switching signal $\sigma \in \Sigma_+ $ and any admissible nonlinearities $\psi$ satisfying \eqref{sectork}. For example,  the nonlinear function $\psi(x)= (\psi_1(x_1)\ \psi_2(x_2))^{\top}$ defined as  $\psi_1(x_1) = 2x_1+ \dfrac{x_1\cos^2 x_1}{1+\sin^2 x_1},\ \psi_2(x_2)= x_2+x_2e^{-x_2^2}$ is obviously admissible and satisfies the sector constraint \eqref{sectork} with $ \beta_i=3, i=1,2$. Then if we take any initial function $\varphi\in \mathcal{C}$ and any switching signal $\sigma \in \Sigma_+ $ the corresponding solution of the above nonlinear system, calculated by the MATLAB code {\bf dde23} decays exponentially to zero. For instance, Fig. 3 shows a signal $\sigma \in \Sigma_+ $ while the corresponding trajectory with the initial condition  $\varphi(\theta)= (\sin(\theta)\ \ \cos(\theta))^{\top}, \ \theta\in [-1,0]$ is shown in Fig. 4. 

\begin{figure}[ht]
		\begin{center}
				\includegraphics[height=6cm,width=11cm]{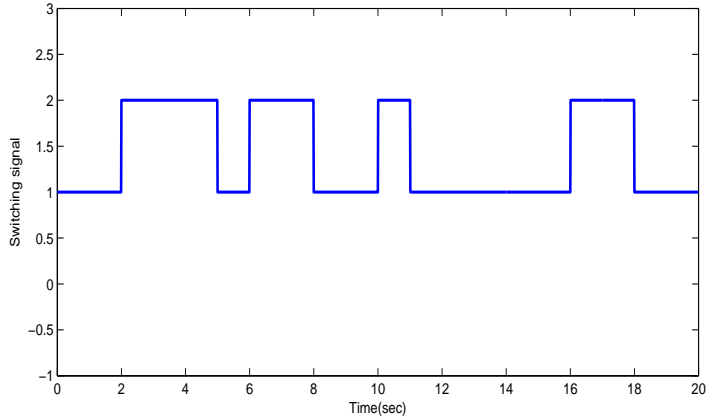}
		\end{center}\vspace*{-0.5cm}
				\caption{\textit{The switching signal $\sigma \in \Sigma_+ $}}
\end{figure}
\vspace{-0.5cm}
\begin{figure}[ht]
		\begin{center}
				\includegraphics[height=6cm,width=11cm]{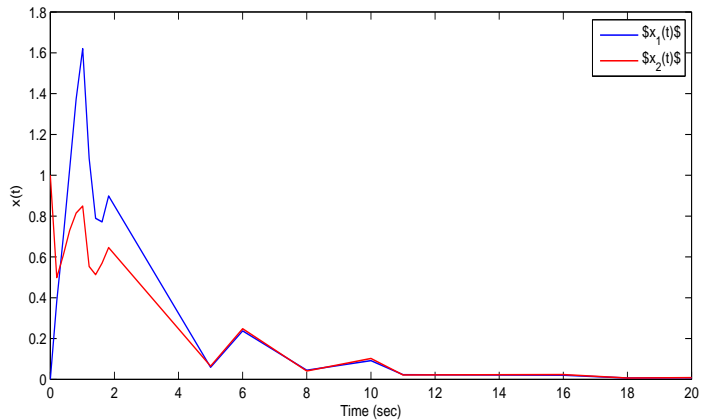}
		\end{center}\vspace{-0.5cm}
				\caption{\textit{The  solution trajectory of Example 2  under any switching signal $\sigma \in \Sigma_+ $}}
\end{figure}
}	
\end{example}
\vspace{-0.5cm}
\section{Conclusion}
In this  paper, we have investigated the exponential stability for nonlinear  switched systems described by time-varying functional differential equations. Some sufficient conditions for exponential stability of the  zero solution have been established, by using the comparison principle and average dwell time switchings. The obtained results apply to different classes of switched time-varying linear and nonlinear systems with delays, yielding verifiable sufficient conditions for exponential stability under arbitrary switching signals with ADT. A special attention is devoted to the absolute stability of switched time-varying systems and a new sufficient condition  in terms of linear inequalities are established to guarantee that
 each solution of the system converges exponentially to zero for any switching signal and  any admissible sector nonlinearities. Even  when applying to the case of time-invariant switched systems our results are reduced to less restrictive conditions for stability, comparing with many previous known results. 
An interesting direction of future research is an extension of
the obtained results to switched discrete-time systems and monotone systems.

\section*{Acknowledgments}
\vskip-0.3cm
 \noindent  This work was partly supported by VAST (Vietnam Academy of  Science and Technology) by the project DLTE001/21-22.

\end{document}